\documentclass{amsart}
\usepackage{color, graphicx, enumerate, amssymb}
\usepackage{hyperref}
\usepackage{epsfig,wrapfig}
\usepackage{pxfonts}
\usepackage{graphicx}

\usepackage{eucal}

\usepackage[utf8]{inputenc}
\usepackage{ulem}

\usepackage{amsmath}

\usepackage{geometry}
\usepackage{enumitem}
\usepackage{amsmath,amssymb,amsthm,upref,graphicx,mathrsfs, enumerate, setspace,xcolor}
\usepackage{esint}
\usepackage{textcomp}

\newtheorem{thm}{Theorem}[section]
\newtheorem{prop}[thm]{Proposition}
\newtheorem{lemma}[thm]{Lemma}
\newtheorem{cor}[thm]{Corollary}
\theoremstyle{definition}
\newtheorem{defi}[thm]{Definition}

\newtheorem{rem}[thm]{Remark}
\newtheorem{eg}[thm]{Example}

\newcommand{\R}{\mathbb{R}}

\newcommand{\eps}{\epsilon}

\newcommand{\del}{\delta}

\newcommand{\fl}{\rightarrow}

\newcommand{\ti}{\times}

\newcommand{\om}{\omega}

\newcommand{\sig}{\sigma}

\newcommand{\al}{\alpha}

\newcommand{\cd}{\cdots}
\newcommand{\p}{\partial}
\newcommand{\n}{\nabla}

\newcommand{\se}{\subseteq}
\newcommand{\fa}{\forall}
\newcommand{\iy}{\infty}

\newcommand{\eqn}[2]{\begin{equation}#2\label{#1}\end{equation}}

\title[]{Weak Form of Differential Equations and Differential Identities}

	\authors
	\author[]{Seyed Ebrahim Akrami}
	\address{ Mathematics Department, Semnan University, IRAN and  Mathematics and Physics Departments, Institute for Research in Fundamental Sciences IPM, IRAN}
	
	\keywords{partial differential equations, weak solution, quantum mechanics, Schr\"{o}dinger equation, Wasserstein space, Stokes's theorem, Euler-Lagrange Equation\\E-mail: akramisa@semnan.ac.ir, akramisa@ipm.ir}
\date{\today}

\begin{document}
\begin{abstract}
Inspired  by quantum mechanics, we introduce a  weak form of solutions  for differential equations and differential identities like Stokes theorem and Euler-Lagrange equation. We show that Schr\"{o}dinger equation is a weak from of the classical Euler-Lagrange equation.
\end{abstract}
	
\maketitle

\section{Introduction}
For partial differential equations various forms of weak solutions are defined, \cite{E}. In this paper, we define a new form of weak solutions for differential equations and differential identities. Motivation comes from quantum mechanics where Newton equation is replaced with the Schr\"{o}dinger equation. In classical mechanics, Newton equation for a particle is 
$$m\ddot{x}(t)+\n U(x(t))=0$$where $m$ is a positive constant, $x(t)$ is a curve in $\R^3$ and $U:\R^3\fl\R$. In quantum mechanics this equation  is replaced with the Schr\"{o}dinger equation 
\eqn{Sch}{i\hbar\frac{\p\psi}{\p t}=-\frac{\hbar^2}{2m}\n\psi+U\psi}where $\psi(t,x)\in\mathbb{C},t\in\R,x\in\R^3.$ There is a big radical change in passing from Newton equation to Schr\"{o}dinger equation. We try to bring these two different equations close to each other. To do this task, we first rewrite Newton equation  as follows
\eqn{xt}{\dot{x}(t)=v(t),}and
\eqn{mvt}{m\dot{v}(t)+\n U(x(t))=0.} The first equation just means that the variable $v(t)$ is the derivative of $x(t)$. Next, using the polar decomposing of the complex function
$\psi(t,x)=R(t,x)e^{\frac{i}{\hbar}S(t,x)}$ and setting $\rho(t,x):=R^2(t,x),V(t,x):=\frac{1}{m}\n S(t,x)$ 
we rewrite the Schr\"{o}dinger equation as follows
\eqn{prhopt}{\frac{\p\rho}{\p t}+\n.(\rho V)=0} and
$$m\Big(\frac{\p V}{\p t}+(V\cdot\n)V\Big)+\n U+\n Q=0$$where the term $Q$ is called quantum potential and it is known that $\int\rho\n Qdx=0$, \cite{H}. Thus 
\eqn{mpVpt1}{\int\rho \Big(m(\frac{\p V}{\p t}+(V\cdot\n)V)+\n U\Big)dx=0.}We call the equations (\ref{prhopt}) and (\ref{mpVpt1}) as \emph{weak form} of the equations (\ref{xt}) and (\ref{mvt}).
This observation suggests to replace a curve $x(t)\in\R^3$ and its derivative $v(t)=dx/ dt$ by probability density $\rho(t,x)$ and vector field $V(t,x)$. We will call the pair $(x(t),v(t))$ as a \emph{strong differentiable function} and the the pair $(\rho(t,x),V(t,x))$ as a \emph{weak differentiable function}. Our goal in this paper is to look for weak form of the well-known differential equations and identities and then solve or prove them. 

Ambrosio et al., \cite{AGS}, have studied differential calculus over a metric measure space and in particular over Wasserstein space of probability measures.  We briefly recall this theory. First, let $(S,d)$ be a complete metric
space. For any curve $\rho:(a, b)\se\R\fl S$ the limit
\eqn{}{|\rho'|(t) := \lim_{h\fl0} \frac{d(\rho(t + h),
		\rho(t))}{|h|}} if exists, is called the \textbf{metric
	derivative} of $\rho$. For any absolutely continuous curve $\rho:(a, b)\se\R\fl S$ the
metric derivative exists for Lebesgue-a.e. $t\in(a, b)$ and we have
$d(\rho(s), \rho(t))\le\int_s^t |\rho'|(r) dr$ for any interval
$(s, t)\se(a, b).$ In the case where $S$ is a Banach space, a curve $\rho:(a,b)\fl S$ is absolutely continuous
if and only it is differentiable in the ordinary sense for
Lebesgue-a.e. $t\in(a, b)$ and we have $\|\rho'(t)\|=|\rho'|(t)$
for Lebesgue-a.e. $t\in(a, b)$. Next in \cite{AGS}, this fact
has been applied to the metric space of probability measures on a
Hilbert space $X$. This space is equipped with the Wasserstein's metric. It is
shown   that the class of absolutely continuous curves
$\rho_t$ in this metric space coincides with solutions of the
continuity equation of physicists. More precisely, given an absolutely
continuous curve $\rho_t$, one can find a Borel time-dependent
velocity field $V_t : X\fl X$ such that $\|V_t\|_{L^p(\rho_t)}\le
|\rho'|(t)$ for a.e. $t$ and the continuity equation
$\frac{\p\rho}{\p t}+\n.(\rho V)=0$ holds. Conversely, if $\rho_t$
solves the continuity equation for some Borel velocity field $V_t$
with $\int_a^b \|V_t\|_{L^p(\rho_t)}dt < \iy$, then $\rho_t$ is an
absolutely continuous curve and $\|V_t\|_{L^p(\rho_t)}\le
|\rho'|(t)$ for a.e. $t \in(a, b).$ As a consequence, we see that
among all velocity fields $V_t$ which produce the same flow
$\rho_t$, there is a unique optimal one with the smallest
$L^p(\rho_t,X)$-norm, equal to the metric derivative of $\rho_t$.
One can  view this optimal field as the ``tangent" vector field to
the curve $\rho_t$. Next, Villani,  Lott, and Sturm studied the concept of Ricci curvature for a metric measure space \cite{LV, V}. The references  \cite{AG, GNT, L1, L2, LV} studied   Riemannian geometry and dynamics on the  Wasserstein space by putting the structure of an informal manifold, first introduced by Otto \cite{O}. 

In this paper, we will work with the  informal differentiable structure on the Wasserstein space. In Section 2, we will define the concept of a weak differentiable function. In Section 3, we will state and prove the weak form of the  Stokes  theorem. In Section 4, we will state the weak form of the Euler-Lagrange equation.

\section{Weakly Differentiable Functions}
This paper is based on the following observation which redefines the Newton-Leibniz derivative of a function.
\begin{thm}\label{1}
A curve $x(t)\in M$ over a manifold $M$ is differentiable  if and only if there exists a vector field $v(t)\in T_{x(t)}M$ along $x(t)$ such that
	\eqn{}{\frac{d}{dt}f(x(t))=df_{x(t)}(v(t))} for all differentiable functions $f$ on $M$ and moreover we have $v(t)=dx(t)/dt$. 
\end{thm}
In this paper, we call the function $x(t)\in M$ a \emph{strong differentiable curve}.
\begin{defi}
A  \textbf{weak curve} over an oriented Riemannian manifold $M$ is a curve $\rho_t$ of  probability densities over $M$, i.e.  $\rho(t,x)\ge0,a\le t\le b,x\in M$,  $\int_M\rho(t,x)dx=1$  for all $a\le t\le b$, and $\rho$ vanishes at infinity with respect to the $x$-variable, where integral is with respect to the Riemannian volume form of $M$. 
\end{defi}
Next, we define the \emph{derivative of a weak curve} by imitating Theorem \ref{1}. 
\begin{defi}
A weak curve $\rho$ over $M$ is called differentiable if there exists a differentiable time-dependent  vector field $V_t$ over $M$ such that
\eqn{fracddtintrhof}{\frac{d}{dt}\int \rho(t,x) f(x)dx=\int \rho(t,x)df(V(t,x))} for all $f\in C_c^\iy(M)$, in the space of all smooth compactly-supported functions on $M$.
\end{defi}
\begin{thm}
A weak curve $\rho$ is differentiable if and only if there exists a differentiable time-dependent  vector field $V_t$ such that
\eqn{f}{\frac{\p \rho}{\p t}+\n.(\rho V)=0.}   
\end{thm}
\begin{proof}
We have $\frac{d}{dt}\int \rho f=\int\frac{\p \rho}{\p t}f$. On the other hand by the divergence theorem and that $f$ has compact support and $\rho$ vanishes at infinity, we have
$$\int \rho df(V)=\int\n.(\rho f V)-\int f\n.(\rho V)=-\int f\n.(\rho V).$$Thus (\ref{fracddtintrhof}) holds if and only if $\int\Big(\frac{\p \rho}{\p t}+\n.(\rho V)\Big)f=0, \fa f$. The later is equivalent with (\ref{f}).
\end{proof}  
	
Note that the vector field $V$ is not unique. Thus   we redefine our concept of differentiability as follows.
\begin{defi}\label{differentiableweakcurve}
A \textbf{differentiable weak curve} is a pair $(\rho_t,V_t)$ where $\rho_t$ is a curve of  probability densities  and $V_t$ is a vector field over $M$ such that  (\ref{f}) holds.
\end{defi}
Now, based on the Introduction, we introduce the central concept of this paper. 
\begin{defi}
By \emph{weakening} a first order differential equation or identity
\eqn{}{G(x(t),\dot{x}(t))=0,} we mean the following equation or identity
\eqn{}{\int_M\rho(t,x)G(x,V(t,x))dx=0} whose unknown or argument is a differentiable weak curve $(\rho,V)$.
\end{defi}

Next, based on the motivation  mentioned in the Introduction, we extend the process of weakening differential equations to higher order equations. 
\begin{defi}
The weak form of the second derivative $\Ddot{x}(t)$ of a strong curve over $M$ for a weak curve $(\rho,V)$ is the following vector field
\eqn{}{\frac{\p V}{\p t}+(V.\n)V.}
\end{defi}
\begin{defi}
By \emph{weakening} a second order differential equation or identity
\eqn{}{G(x(t),\dot{x}(t),\Ddot{x}(t))=0,} we mean the following equation or identity
\eqn{}{\int_M\rho G(x,V,\frac{\p V}{\p t}+(V.\n)V)dx=0} whose unknown or argument is a differentiable weak function $(\rho,V)$.
\end{defi}
 
Next, we define a several variable differentiable weak function.
\begin{defi}\label{multi differentiable np }
	A \emph{differentiable weak   function} from an open subset
	$U\se\R^m$ into  $M$ is a system
	$(\rho,V_1,\cd,V_m)$ including a differentiable function $\rho:U\ti M\fl[0,\iy)$, vanishing at infinity with respect to $M$,  together
	with differentiable vector fields $V_j:U\ti M\fl TM$, $1\le
	j\le m,$ such that \eqn{}{\frac{\p \rho}{\p u_j}+\n.(\rho V_j)=0.}
\end{defi}
A weak function $\rho:U\ti M\fl[0,\iy)$ is a replacement for a strong
function $f:U\fl M$ and the vector fields $V_i$ are replacements
for the partial derivatives $V_i=\frac{\p f}{\p u_i}$. For the strong function $f$ we have
$\frac{\p^2 f}{\p u_i\p u_j}=\frac{\p^2 f}{\p u_j\p u_i}$. Namely, for  $V_i=\frac{\p f}{\p u_i}$ we have  \eqn{pAi}{\frac{\p V_i}{\p
		u_j}=\frac{\p V_j}{\p u_i}.}
For weak functions, we have the
following result.
\begin{thm}If $(\rho,V_1,\cd,V_m)$ is  a differentiable weak  function from an open region $U\se\R^m$ to $M$ then
	\eqn{}{\rho\Big(\frac{\p V_i}{\p u_j}-\frac{\p V_j}{\p
			u_i}-[V_i,V_j]\Big)=0.}Here $[V_i,V_j]$ means the Lie bracket of
	vector fields $V_i$ and $V_j$ with respect to $M$.
\end{thm}
\begin{proof}For any compactly-supported function $g$ on $M$, since $\rho$ vanishes at infinity, by the divergence theorem we have
	\begin{eqnarray}\frac{\p}{\p u_i}\int \rho g&=&\int \frac{\p \rho }{\p u_i}g\nonumber\\
		&=&-\int g\n.(\rho V_i)\nonumber\\
		&=&-\int\n(\rho gV_i)+\int \rho V_ig\nonumber\\
		&=&\int \rho V_ig.\nonumber\end{eqnarray}Thus
	\begin{eqnarray}\frac{\p^2}{\p u_j\p u_i}\int \rho g&=&\int\frac{\p \rho }{\p u_j}V_ig+\int \rho \frac{\p V_i}{\p u_j}g\nonumber\\
		&=&-\int \n.(\rho V_j)V_ig+\int \rho \frac{\p V_i}{\p u_j}g\nonumber\\
		&=&-\int \n.(\rho (V_ig)V_j)+\int \rho V_jV_ig+\int \rho \frac{\p V_i}{\p u_j}g\nonumber\\
		&=&\int \rho (\frac{\p V_i}{\p u_j}+V_jV_i)g.\nonumber\end{eqnarray}
	Hence 
	\begin{eqnarray}0&=&\frac{\p^2}{\p u_j\p u_i}\int \rho g-\frac{\p^2}{\p u_i\p u_j}\int \rho g\nonumber\\&=&\int \rho (\frac{\p V_i}{\p u_j}+V_jV_i-\frac{\p V_j}{\p u_i}-V_iV_j)g\nonumber\\&=&\int \rho (\frac{\p V_i}{\p u_j}-\frac{\p V_j}{\p u_i}+[V_i,V_j])g.\nonumber \end{eqnarray} Since $g$ is arbitrary, the desired identity is obtained.
\end{proof}
\begin{rem}
Without using the divergence theorem, we can prove at most the following identity
\eqn{}{\n.\Big(\rho(\frac{\p V_i}{\p u_j}-\frac{\p V_j}{\p
		u_i}-[V_i,V_j])\Big)=0.}To prove the later, first, we state an important identity that we discovered recently and we have not seen before elsewhere; one can show by some calculation  that for any two vector fields $V$ and $W$ and any function $f$ on $M$ we have
\eqn{}{\n.\Big(\n.(f W)V\Big)-\n.\Big(\n.(f V)W\Big)=\n.(f[W,V]).}  Now using this identity  we have	
\begin{eqnarray}
\frac{\p^2\rho}{\p u_j\p u_i}&=&-\frac{\p}{\p u_j}\n.(\rho V_i)\nonumber\\&=&-\n.(\frac{\p\rho}{\p u_j}V_i)-\n.(\rho\frac{\p V_i}{\p u_j})\nonumber\\&=&\n.\Big(\n.(\rho V_j)V_i\Big)-\n.(\rho\frac{\p V_i}{\p u_j}).\nonumber
\end{eqnarray}
Thus 
\begin{eqnarray}
0&=&\frac{\p^2\rho}{\p u_j\p u_i}-\frac{\p^2\rho}{\p u_i\p u_j}\nonumber\\&=&\n.\Big(\n.(\rho V_j)V_i\Big)-\n.\Big(\n.(\rho V_i)V_j\Big)-\n.(\rho\frac{\p V_i}{\p u_j})+\n.(\rho\frac{\p V_j}{\p u_i})\nonumber\\&=&-\n.(\rho[V_i,V_j])-\n.(\rho\frac{\p V_i}{\p u_j})+\n.(\rho\frac{\p V_j}{\p u_i})\nonumber\\&=&\n.\Big(\rho(\frac{\p V_j}{\p u_i}-\frac{\p V_i}{\p u_j}-[V_i,V_j])\Big).\nonumber
\end{eqnarray}
	
\end{rem}

\begin{defi}\label{multi differentiable np }
	A \emph{differentiable weak   function} from an open subset
	$U\se\R^m$ into  $M$ is a system
	$(\rho,V_1,\cd,V_m)$ including a differentiable function $\rho:U\ti M\fl[0,\iy)$   together
	with differentiable vector fields $V_j:U\ti M\fl TM$, $1\le
	j\le m,$ such that 
	\eqn{}{\frac{\p \rho}{\p u_j}+\n.(\rho V_j)=0.} We sometimes use the brief notation $(\rho,V):U\fl M$ to denote a weak differentiable function from $U$ to $M$. 
\end{defi}

\begin{eg}
	Let $A\in M_{n\ti m}(\R)$ be a matrix of size $n\ti m$ with real entries whose columns are  $A_i,1\le i\le m,$ and let $\sig$ be a probability density over  $\R^n$. Then the pair $(\rho,V_1,\cd,V_m)$ given below is a differentiable weak function from $\R^m$ into $\R^n$.
	\eqn{}{\rho(x,y):=\sig(y-Ax),~~~~~~V_i(x,y):=A_i.}
\end{eg}
\begin{proof}We have $\frac{\p\rho}{\p x_i}=-\n \sig(y-Ax).A_i$ and $\n^y.(\rho(x,y)V_i(x,y))=\n^y.(\sig(y-Ax)A_i)=\n\sig(y-Ax).A_i$. 
\end{proof}

\begin{defi}\label{smoouncerfun}
	A \emph{differentiable weak  function} from a manifold $P$
	into an oriented Riemannian manifold $M$ is a pair $(\rho,V)$
	including a differentiable   function $\rho:P\ti M\fl[0,\iy)$, together with a differentiable map
	\eqn{}{\boxed{V:TP\ti M\fl TM},~~~~~~V(w,q)\in T_qM} $\fa w\in TP,q\in M$,
	such that $V$ is linear concerning the first variable and for any coordinate system
	$u=(u_1,\cd,u_m)\in U\se\R^m$ for $P$ \eqn{pphi}{\frac{\p \rho}{\p
			u_j}+\n.(\rho V_j)=0} where $V_j:U\ti M\fl
	TM,V_j(u,q):=V(\frac{\p}{\p u_j}(u),q)\in T_qM.$ Here the
	divergence operator $\n.$ is relative to the variable $q\in M$. 
	
\end{defi}
\begin{prop}\label{dprho}
	The  equations (\ref{pphi}) does not depend on the coordinate systems of $P$.
\end{prop}
\begin{proof}Let $W_i(v):=V_{p}(\frac{\p}{\p v_i}(p))$ where $(v_1,\cd, v_m)$ is another coordinate system for $P$.  Then since $\frac{\p}{\p
		v_i}=\sum_j\frac{\p u_j}{\p v_i}\frac{\p}{\p u_j}$ we get
	$W_i(v)=\sum_j\frac{\p u_j}{\p v_i}V_j(u)$. Thus
	\begin{eqnarray}& &\frac{\p\rho(v,q)}{\p v_i}+\n.(\rho(v,q)W_i(v,q))
		\nonumber\\&=& \sum_j\frac{\p u_j}{\p v_i}\frac{\p\rho(u,q)}{\p
			u_j}+\n.(\rho(u,q)\sum_j\frac{\p u_j}{\p
			v_i}V_j(u,q))\nonumber\\&=&\sum_j\frac{\p u_j}{\p
			v_i}[\frac{\p\rho(u,q)}{\p u_j}+\n.(\rho(u,q)V_j(u,q))].\nonumber
	\end{eqnarray}Thus since the matrix $(\frac{\p u_j}{\p
		v_i})$ is invertible we conclude that $\frac{\p\rho(v,q)}{\p
		v_i}+\n.(\rho(v,q)W_i(v,q))=0$ if and only if
	$\frac{\p\rho(u,q)}{\p u_i}+\n.(\rho(u,q)V_i(u,q))=0.$
	
\end{proof}

\section{Weak Form of Stokes  Theorem}
In the strong  case when we have a differentiable map $F:Q\fl M$ from a $k$-dimensional oriented manifold $Q$ to $M$ then we can pull back any $k$-differential form $\om$ on $M$ to a $k$-dimensional differential form $F^*\om$ on $Q$ and then integrate over $Q$ as $\int_QF^*\om$. Then we have the Stokes  theorem  on map $F$
\eqn{}{\int_QF^*d\om=\int_{\p Q}F^*\om.}The proof is easy: $\int_QF^*d\om=\int_QdF^*\om=\int_{\p Q}F^*\om$ by the Stokes  theorem on $Q$. Now we are going to prove a weak version of this theorem. 
\begin{defi}\label{npcpull}
	Let $F=(\rho,V_1,\cd,V_k)$ be a weak map from $Q$ to $M$. The pullback of $k$-differential form $\om$ on $M$ to a $k$-differential form $F^*\om$ on $Q$ is defined by
	\eqn{}{(F^*\om)(q;w_1,\cd,w_k):=\int_M\rho(q,p)\om(p;V(w_1)(p),\cd,V(w_k)(p))} for $q\in Q,p\in M,w_i\in T_qQ.$
\end{defi}
\begin{thm}
	\eqn{}{F^*d\om=dF^*\om.}
\end{thm}

\begin{proof}To show $F^*d\om=dF^*\om$, it is enough to apply both sides  to $(\frac{\p}{\p u_{l_1}},\cd,\frac{\p}{\p u_{l_k}})$ where $(u_1,\cd,u_k)$ is a coordinate system for $Q$. As before we set $V_i=V(\frac{\p}{\p u_i}).$
\begin{eqnarray}
		& &(F^*d\om)(q;\frac{\p}{\p u_{l_1}},\cd,\frac{\p}{\p u_{l_k}})\nonumber\\&=&\int_M\rho(q)d\om(V(\frac{\p}{\p u_{l_1}}),\cd,V(\frac{\p}{\p u_{l_k}}))\nonumber\\&=&\int_M\rho(q)d\om(V_{l_1},\cd,V_{l_k})\nonumber\\&=&\sum_{i=1}^{k}(-1)^{i-1}\int_M\rho(q)V_{l_i}\om(V_{l_1},\cd,\widehat{V_{l_i}},\cd,V_{l_k})\nonumber\\&+&\sum_{i< j}(-1)^{i+j}\int_M\rho(q)\om([V_{l_i},V_{l_j}],V_{l_1},\cd,\widehat{V_{l_i}},\cd,\widehat{V_{l_j}},\cd V_{l_k})\nonumber\\&=&
		-\sum_{i=1}^{k}(-1)^{i-1}\int_M\n.(\rho V_{l_i})\om(V_{l_1},\cd,\widehat{V_{l_i}},\cd,V_{l_k})\nonumber\\&+&\sum_{i< j}(-1)^{i+j}\int_M\rho(q)\om([V_{l_i},V_{l_j}],V_{l_1},\cd,\widehat{V_{l_i}},\cd,\widehat{V_{l_j}},\cd V_{l_k})\nonumber\\&=&
		\sum_{i=1}^{k}(-1)^{i-1}\int_M\frac{\p\rho}{\p u_{l_i}}\om(V_{l_1},\cd,\widehat{V_{l_i}},\cd,V_{l_k})\nonumber\\&+&\sum_{i< j}(-1)^{i+j}\int_M\rho(q)\om([V_{l_i},V_{l_j}],V_{l_1},\cd,\widehat{V_{l_i}},\cd,\widehat{V_{l_j}},\cd V_{l_k})\nonumber\\&=&
		\sum_{i=1}^{k}(-1)^{i-1}\Big(\frac{\p}{\p u_{l_i}}\int_M\rho\om(V_{l_1},\cd,\widehat{V_{l_i}},\cd,V_{l_k})\nonumber\\&-&
		\int_M\rho\frac{\p}{\p u_{l_i}}(\om(V_{l_1},\cd,\widehat{V_{l_i}},\cd,V_{l_k}))\Big)
		\nonumber\\&+&\sum_{ i< j}(-1)^{i+j}\int_M\rho(q)\om([V_{l_i},V_{l_j}],V_{l_1},\cd,\widehat{V_{l_i}},\cd,\widehat{V_{l_j}},\cd V_{l_k})\nonumber\\&=&
		\sum_{i=1}^{k}(-1)^{i-1}\Big(\frac{\p}{\p u_{l_i}}\int_M\rho\om(V_{l_1},\cd,\widehat{V_{l_i}},\cd,V_{l_k})\nonumber\\&-&
		\sum_{ j<i}(-1)^{j-1}\int_M\rho\om(\frac{\p V_{l_j}}{\p u_{l_i}},V_{l_1},\cd,\widehat{V_{l_j}},\cd,\widehat{V_{l_i}},\cd,V_{l_k})\nonumber\\&-&\sum_{ i<j}(-1)^{j-2}\int_M\rho\om(\frac{\p V_{l_j}}{\p u_{l_i}},V_{l_1},\cd,\widehat{V_{l_i}},\cd,\widehat{V_{l_j}},\cd,V_{l_k})\Big)
		\nonumber\\&+&\sum_{i< j }(-1)^{i+j}\int_M\rho(q)\om([V_{l_i},V_{l_j}],V_{l_1},\cd,\widehat{V_{l_i}},\cd,\widehat{V_{l_j}},\cd V_{l_k})\nonumber\\
		&=&\sum_{i=1}^{k}(-1)^{i-1}\frac{\p}{\p u_{l_i}}\int_M\rho\om(V_{l_1},\cd,\widehat{V_{l_i}},\cd,V_{l_k})\nonumber\\&+&
		\sum_{i<j}(-1)^{i+j}\int_M\rho\om(\frac{\p V_{l_j}}{\p u_{l_i}}-\frac{\p V_{l_i}}{\p u_{l_j}},V_{l_1},\cd,\widehat{V_{l_j}},\cd,\widehat{V_{l_i}},\cd,V_{l_k})
		\nonumber\\&+&\sum_{i< j }(-1)^{i+j}\int_M\rho(q)\om([V_{l_i},V_{l_j}],V_{l_1},\cd,\widehat{V_{l_i}},\cd,\widehat{V_{l_j}},\cd V_{l_k})\nonumber
		\end{eqnarray}
		\begin{eqnarray}&=&
		\sum_{i=1}^{k}(-1)^{i-1}\frac{\p}{\p u_{l_i}}\int_M\rho\om(V_{l_1},\cd,\widehat{V_{l_i}},\cd,V_{l_k})\nonumber\\&+&
		\sum_{i<j}(-1)^{i+j}\int_M\rho\om(\frac{\p V_{l_j}}{\p u_{l_i}}-\frac{\p V_{l_i}}{\p u_{l_j}}+[V_{l_i},V_{l_j}],V_{l_1},\cd,\widehat{V_{l_j}},\cd,\widehat{V_{l_i}},\cd,V_{l_k})\nonumber\\&=&
		\sum_{i=1}^{k}(-1)^{i-1}\frac{\p}{\p u_{l_i}}\int_M\rho\om(V_{l_1},\cd,\widehat{V_{l_i}},\cd,V_{l_k})\nonumber\\&=&
		\sum_{i=1}^{k}(-1)^{i-1}\frac{\p}{\p u_{l_i}}F^*\om(\frac{\p}{\p u_{l_1}},\cd,\widehat{\frac{\p}{\p u_{l_i}}},\cd,\frac{\p}{\p u_{l_k}})\nonumber\\&+&\sum_{i<j}(-1)^{i+j}F^*\om([\frac{\p}{\p u_{l_i}},\frac{\p}{\p u_{l_j}}],\frac{\p}{\p u_{l_1}},\cd,\widehat{\frac{\p}{\p u_{l_i}}},\cd,\widehat{\frac{\p}{\p u_{l_j}}},\cd \frac{\p}{\p u_{l_k}})\nonumber\\&=&dF^*\om(q;\frac{\p}{\p u_{l_1}},\cd,\frac{\p}{\p u_{l_k}}).\nonumber
	\end{eqnarray}
\end{proof} 
\begin{thm}\label{npcstokes}(\emph{Weak Stokes  theorem}) For any weak map, $F=(\rho,V)$ from a $k$-dimensional manifold $Q$ to a manifold $M$ and any $k$-form $\om$ over $M$
	\eqn{}{\int_QF^*d\om=\int_{\p Q}F^*\om.}
\end{thm}
\begin{proof}
	$\int_QF^*d\om=\int_QdF^*\om=\int_{\p Q}F^*\om.$
\end{proof}

Let us see what this theorem says in Euclidean spaces. The ordinary (strong) Stokes  theorem for a  surface $S$ in $\R^3$ parameterized
by an ordinary (strong) map $X:D\subseteq \R^2\fl \R^3$ and for a vector field $F(x,y,z)$
in $\R^3$ is given by  
\eqn{}{\int_D \Big((\n\ti F)\circ X\Big).(U\ti V)dq=\int_{\partial
		D}(F\circ X).(Udu+Vdv)}
where  $U=\frac{\p X}{\p u}, V=\frac{\p X}{\p v}$ and $dq=dudv$. Next, consider a weak  parameterized surface, i.e. a weak differentiable function $(\rho,U,V)$ form a domain $D\subseteq\R^2$ into $\R^3$, including  a positive function $\rho(q;p)\ge0$, vanishing at infinity concerning the variable $p$, and two vector fields $U=U(q;p)\in\R^3$,$V=V(q;p)\in\R^3$,
$q=(u,v)\in D\subseteq\R^2, p=(x,y,z)\in\R^3$ such that   $$\frac{\partial\rho}{\partial
	u}+\nabla.(\rho U)=0,~~~\frac{\partial\rho}{\partial
	v}+\nabla.(\rho V)=0.$$ 
Then we have
\eqn{}{ \int_{\R^3}\int_D\rho (\n\ti F).(U\ti V)dqdp=\int_{\R^3}\int_{\partial
		D}\rho F.(Udu+Vdv)dp}
where $dq=dudv$ and $dp=dxdydz$.

\section{Weak Form of Euler-Lagrange Equation}
Classical mechanics over   Wasserstein space of probability measures over a manifold has been studied extensively, \cite{AG,GNT,GS,AG2,Ca}. In this section, we are going to give another viewpoint on this subject. In classical mechanics over strong states, i.e. points over a finite-dimensional configuration manifold $M$, we study the following action functional
\eqn{Sal}{S(\al)=\int_a^bL(\al(t),\dot{\al}(t))dt} over all paths $\al(t)\in M,a\le t\le b,$ and $L:TM\fl\R$ is a Lagrangian. The Lagrangian can have an arbitrary, i.e. non-restrictive, dependence to both $\al$ and its derivative. Now, in this paper we are going to replace the space $M$ with the space $W^\iy(M)$ of smooth probability measures over $M$, i.e. smooth Wasserstein space. 
The first natural question would be what is the correct counterpart of the above action functional in the case of the space $W^\iy(M)$? The Lagrangian in the above action is a function of $\al$ and its derivative. So the Lagrangian over $W^\iy(M)$ must be a function of a path $d\mu_t=\rho(t)dx\in W^\iy(M)$ and its derivative. In this paper we work with probability measures which have density with respect to the volume form of Riemannian manifold $M$ where its Riemannian volume form is denoted by $dx$. We know that the time-derivative of $\rho(t)$ is a vector field $V(t)$ over $M$ satisfying the continuity equation in the sense of metric measure space theory. So the Lagrangian must be  a function of $\rho$ and $V$. The first candidate for an action functional  over $W^\iy(M)$ 
is 
\eqn{}{S(\rho,V)=\int_a^b\int_M\rho(t,x)L(x,V(t,x))dxdt,}where $L:TM\fl\R$ is a classical Lagrangian over $M$  and $dx$ is the Riemannian volume form of $M$. But this is too restrictive Lagrangian since its dependence to the density $\rho$ is very simple unlike the action (\ref{Sal}) which has nonrestrictive dependence to $\al$. Thus we must add some extra term to the Lagrangian which depends to $\rho$ and perhaps to its spatial derivatives.

In \cite{GNT} the action functional is of the form
\eqn{}{S(\mu,V)=\frac{1}{2}\int_a^b\int_{\R^n}\Arrowvert V(t,x)\Arrowvert^2 d\mu_t(x)dt-\int_a^bU(\mu_t)dt} where $\mu_t$ is a path of probability measures and $V_t$ is vector field over $\R^n$. Here, we do not assume that the continuity equation between $\mu_t$ and $V_t$ holds.

In \cite{Ca} the action functional is of the following form
\eqn{SmuVFG}{S(\mu,V)=\int_a^b\int_{\R^n}L(x,V(t,x))d\mu_t(x)dt+\int_a^bF(\mu_t)dt+G(\mu(b))}where $\mu_t$ is a curve of probability measures over $\R^n,V_t$ is a curve of vector fields over $\R^n, L:T\R^n\fl\R,$ and $F$ as well as $G$ are functionals over the space of probability measures and the continuity equation between $\mu_t$ and $V_t$ holds, i.e. $\frac{\p\mu_t}{\p t}+\n.(\mu_tV_t)=0$ in the sense of distributions.

In this section, we study the following  action functional
\eqn{SrhoV}{S(\rho,V)=\int_a^b\int_M\rho(t,x)\Big(L(x,V(t,x))+F(\rho,\p_i\rho,\p^2_{ij}\rho,\p^3_{ijk}\rho,\cd) \Big)dxdt}
for some function $F(y,y_i,y_{ij}, y_{ijk},\cd)$ and where $\p_i=\frac{\p}{\p x_i}, \p^2_{ij}=\frac{\p}{\p x_i\p x_j},$ etc.
Here, we assume that $(\rho,V)$ is a weak differentiable function.

For simplicity, we will work on space $M=\R^n$ but the result can be extended to an arbitrary manifold. In order to be able to apply the least action principle to the  action \ref{SrhoV}, we first need to define a perturbation of the differentiable curve $(\rho_t,V_t)$. To do this task, we recall the meaning of a perturbation of a strong  differentiable curve $x(t)\in M,a\le t\le b$. It is just a two-variable function $x(t,s)\in M,a\le t\le b,-\eps\le s\le \eps,$ such that $x(t,0)=x(t),x(a,s)=x(a)$ and $x(b,s)=x(b).$ Notice, that we have the following identity
\eqn{ts}{\frac{\p^2x}{\p s\p t}=\frac{\p^2x}{\p t\p s}.}
	
We define a variation of a differentiable weak curve $(\rho_t,V_t),a\le t\le b,$ over $M$  to be a triple $(\rho_{t,s}(x),V_{t,s}(x), W_{t,s}(x))$ for some parameter $s$ in a small interval around zero,  satisfying 
\eqn{rhoVW}{\frac{\p\rho}{\p t}+\n.(\rho V)=0,~~~~~~~~~~\frac{\p\rho}{\p s}+\n.(\rho W)=0,} and moreover for $s=0$ we recover the original pair $(\rho_t,V_t)$ and also $\rho_{a,s}(x)=\rho_{a}(x),\rho_{b,s}(x)=\rho_{b}(x),W_{a,s}(x)=W_{b,s}(x)=0,\fa x$.
The notations in the following lemma will be used in our calculus of variation.
\begin{lemma}
	If $A$ and $B$ be two vectors in $\R^n$ and $X$ is a $n\ti n$ matrix we set $AB$ be the inner product and $XA$ be the vector $XA^T$. Then we have\\
	i) $(XA)B=(XB)A$,\\
	ii) $(AX)B=A(XB).$
\end{lemma}
The proof is obvious.
\begin{thm}\label{wlELE} By the variation of the action functional (\ref{SrhoV}) under the above meaning of variation we arrive to the following equation which we call \emph{weak Euler-Lagrange Equation}
\eqn{WeakELE}{\rho\Big((\frac{\p}{\p t}+V\frac{\p}{\p
		x})\frac{\p L}{\p \dot{x}}-\frac{\p L}{\p
		x}-\frac{\p}{\p x}F-\frac{\p}{\p x}(\rho\frac{\p F}{\p y}-\frac{\p}{\p x_i}(\rho\frac{\p F}{\p y_i})+\frac{\p^2}{\p x_i\p x_j}(\rho\frac{\p F}{\p y_{ij}})-\cd)\Big)=0} where $\n$ is the gradient operator with respect to the  variable $x\in M$ and in general for given vector fields $X$ and function $f$, by $(X.\n)f$ we mean the derivative of function $f$ in the direction of vector field $X$. Here for repeated indexes we use Einstein summation rule. In particular, if the following identity holds for all function $\rho(x)$
\eqn{F}{\rho\frac{\p F}{\p y}-\frac{\p}{\p x_i}(\rho\frac{\p F}{\p y_i})+\frac{\p^2}{\p x_i\p x_j}(\rho\frac{\p F}{\p y_{ij}})-\cd=0} then the weak Euler-Lagrange equation becomes
\eqn{WeakELE1}{\frac{\p\rho}{\p t}+\n.(\rho V)=0,\,\,\,\,\,\rho\Big((\frac{\p }{\p t}+V.\n)\frac{\p L}{\p \dot{x}}(x,V)-\frac{\p L}{\p x}(x,V)-\frac{\p}{\p x}F\Big)=0.}
\end{thm}
\begin{proof}
	We set
	$$S_1:=\int_a^b\int_M\rho
	L(x,V)dxdt$$ and $$S_2:=\int_a^b\int_M\rho F(\rho,\p_i\rho,\p^2_{ij}\rho,\p^3_{ijk}\rho,\cd)
	dxdt.$$
	Let $(\rho_{t, s }(x),V_{t, s }(x), W_{t, s }(x))$ be a variation of the pair $(\rho_t,V_t)$. We set  $K:=L(x,V(t, s ,x))$. Based on the notation introduced in the previous Lemma, we compute the derivative 
\begin{eqnarray} 
& &\frac{dS_1}{d s }=\iint
		\Big(\frac{\p\rho}{\p  s }K+\rho \frac{\p L}{\p \dot{x}} \frac{\p V}{\p
			 s }\Big)dxdt \nonumber\\&=&\iint \rho \Big(\frac{\p K}{\p
			x} W+ \frac{\p L}{\p \dot{x}} \frac{\p V}{\p
			 s }\Big)dxdt\nonumber\\&=&\iint \rho \Big(\frac{\p K}{\p
			x} W  +\frac{\p L}{\p \dot{x}} \frac{\p W}{\p
			t}+ \frac{\p L}{\p \dot{x}}  [V,W]\Big)dxdt\nonumber\\&=&\iint \Big(\rho \frac{\p K}{\p
			x} W+\frac{\p}{\p t}(\rho \frac{\p L}{\p \dot{x}}W)-\rho (\frac{\p}{\p
			t}\frac{\p L}{\p \dot{x}} )W- \frac{\p\rho}{\p t}
		\frac{\p L}{\p \dot{x}}W+ \rho\frac{\p L}{\p\dot{x}}  [V,W]\Big)dxdt\nonumber\\&=&\iint \rho\Big( \frac{\p K}{\p
			x} W- (\frac{\p}{\p
			t}\frac{\p L}{\p \dot{x}} ) W-((\frac{\p}{\p
			x}(\frac{\p L}{\p \dot{x}}  W)) V + \frac{\p L}{\p \dot{x}}  [V,W]\Big)dxdt\nonumber\\&=&\iint \rho\Big( \frac{\p K}{\p
			x} W- (\frac{\p}{\p
			t}\frac{\p L}{\p \dot{x}} ) W-((\frac{\p}{\p x}\frac{\p L}{\p \dot{x}} ) W) V -(\frac{\p L}{\p \dot{x}}  \frac{\p W}{\p x}) V + \frac{\p L}{\p \dot{x}} ( \frac{\p W}{\p x}V-\frac{\p
			V}{\p x}W)\Big)dxdt\nonumber\\&=&\iint \rho \Big(\frac{\p K}{\p
			x} -\frac{\p}{\p
			t}\frac{\p L}{\p \dot{x}}  -(\frac{\p}{\p
			x}\frac{\p L}{\p \dot{x}} )V  -\frac{\p L}{\p \dot{x}} \frac{\p
			V}{\p x}\Big)Wdxdt\nonumber
	\end{eqnarray}
	But
	$\frac{\p K}{\p
		x}=\frac{\p L}{\p x}+\frac{\p L}{\p
		\dot{x}} \frac{\p V}{\p x}$, and thus
	\eqn{}{\del S_1=\frac{d S}{d s }|_{ s =0}
		=\iint \rho\Big(\frac{\p L}{\p
			x}-(\frac{\p}{\p t}+V\frac{\p}{\p
			x})\frac{\p L}{\p \dot{x}}\Big)Wdxdt.} 
	By applying integration by parts and the fact that $\rho$ and thus its derivatives  vanish at infinity, we have
	\begin{eqnarray} & &\frac{dS_2}{d s }=\iint
		\Big(\frac{\p\rho}{\p  s }F+\rho(\frac{\p F}{\p y}\frac{\p\rho}{\p  s }+\frac{\p F}{\p y_i}\frac{\p^2\rho}{\p  s \p x_i}+\frac{\p F}{\p y_{ij}}\frac{\p^3\rho}{\p  s \p x_i\p x_j}+\cd)\Big)dxdt\nonumber\\&=&\iint
		\Big(\frac{\p\rho}{\p  s }F+\rho\frac{\p F}{\p y}\frac{\p\rho}{\p  s }-(\frac{\p}{\p x_i}(\rho\frac{\p F}{\p y_i}))\frac{\p\rho}{\p  s }+(\frac{\p^2}{\p x_i\p x_j}(\rho\frac{\p F}{\p y_{ij}}))\frac{\p\rho}{\p  s }-\cd\Big)dxdt\nonumber\\&=&\iint
		\frac{\p\rho}{\p  s }\Big(F+\rho\frac{\p F}{\p y}-\frac{\p}{\p x_i}(\rho\frac{\p F}{\p y_i})+\frac{\p^2}{\p x_i\p x_j}(\rho\frac{\p F}{\p y_{ij}})-\cd\Big)dxdt\nonumber\\&=&
		\iint
		\rho\frac{\p}{\p x}\Big(F+\rho\frac{\p F}{\p y}-\frac{\p}{\p x_i}(\rho\frac{\p F}{\p y_i})+\frac{\p^2}{\p x_i\p x_j}(\rho\frac{\p F}{\p y_{ij}})-\cd\Big)Wdxdt\nonumber
	\end{eqnarray}
	Thus 
	$$\frac{d S}{d s }|_{ s =0}
	=\iint \rho\Big(\frac{\p L}{\p
		x}-(\frac{\p}{\p t}+V\frac{\p}{\p
		x})\frac{\p L}{\p \dot{x}}+\frac{\p}{\p x}F+\frac{\p}{\p x}(\rho\frac{\p F}{\p y}-\frac{\p}{\p x_i}(\rho\frac{\p F}{\p y_i})+\frac{\p^2}{\p x_i\p x_j}(\rho\frac{\p F}{\p y_{ij}})-\cd)\Big)Wdxdt.$$ 
	Since the variation is arbitrary, we get the desired result. 
\end{proof}
The following fact is due to Iranian theoretical physicists   Mehdi Golshani, Mahdi Atiq and Mozafar Karamian, \cite{AKG} which we have adapted with our theory.
\begin{cor}
When $L(x,v)=\frac{1}{2m}\|v\|^2-U(x)$ and \eqn{}{F(y,y_i,y_{ij})=-\frac{\hbar^2}{2m}(-\frac{\sum_iy_i^2}{4y^2}+\frac{\sum_iy_{ii}}{2y})}  then first of all the identity (\ref{F}) holds for all $\rho$ and secondly if in the weak Euler-Lagrange equation the vector field $V$ is of the form $V=\frac{\n S}{m}$ for some function $S$ then the weak Euler-Lagrange equations convert to the  Schr\"{o}dinger equation (\ref{Sch}), under the transformation
\eqn{psi}{\psi:=\sqrt{\rho}e^{\frac{i}{\hbar}S}.} In this case we have, \eqn{F1}{F(\rho,\p_i\rho,\p_{ij}^2\rho)=-\frac{\hbar^2}{2m}\frac{\n^2\sqrt{\rho}}{\sqrt{\rho}}} 
which is called quantum potential in the literature of Bohmian quantum mechanics.
\end{cor}
\begin{proof}
We have $\frac{\p F}{\p y}=\frac{1}{2}y^{-3}\sum_iy_i^2-\frac{1}{2}y^{-2}\sum_iy_{ii},$ $\frac{\p F}{\p y_i}=\frac{-1}{2}y^{-2}y_i$ and $\frac{\p F}{\p y_{ij}}=\frac{1}{2}\del_{ij}y^{-1}.$  Thus  	$\frac{\p F}{\p y}(\rho,\p_i\rho,\p_{ij}^2\rho)=\frac{1}{2}\rho^{-3}\sum_i(\p_i\rho)^2-\frac{1}{2}\rho^{-2}\sum_i\p_{ii}^2\rho,$ $\frac{\p F}{\p y_i}(\rho,\p_i\rho,\p_{ij}^2\rho)=\frac{-1}{2}\rho^{-2}\p_i\rho$ and $\frac{\p F}{\p y_{ij}}(\rho,\p_i\rho,\p_{ij}^2\rho)=\frac{1}{2}\del_{ij}\rho^{-1}.$ Hence, $\rho \frac{\p F}{\p y}(\rho,\p_i\rho,\p_{ij}^2\rho)-\sum_i\p_i(\rho\frac{\p F}{\p y_i}(\rho,\p_i\rho,\p_{ij}^2\rho))+\sum_{ij}\p_{ij}^2(\rho\frac{\p F}{\p y_{ij}}(\rho,\p_i\rho,\p_{ij}^2\rho))=\frac{1}{2}\rho^{-2}\sum_i(\p_i\rho)^2-\frac{1}{2}\rho^{-1}\sum_i\p_{ii}^2\rho+\frac{1}{2}\sum_i(\p_i\rho^{-1}\p_i\rho)=
\frac{1}{2}\rho^{-2}\sum_i(\p_i\rho)^2-\frac{1}{2}\rho^{-1}\sum_i\p_{ii}^2\rho+\frac{1}{2}(-\rho^{-2}\sum_i(\p_i\rho)^2+\rho^{-1}\sum_i\p_{ii}^2\rho)=0$. Also, one can easily check (\ref{F1}).

Next, the equation (\ref{WeakELE})
in this case becomes
\eqn{}{\rho\Big((\frac{\p}{\p t}+V.\n)V+\n U+\frac{\hbar^2}{2m}\n\frac{\n^2\sqrt{\rho}}{\sqrt{\rho}}\Big)=0.} One can easily check that the later together with the equation $\frac{\p \rho}{\p t}+\n.(\rho V)=0$ when $V=\frac{\n S}{m}$ is equivalent with Schr\"{o}dinger equation under the transformation (\ref{psi}).
\end{proof}

\section*{Acknowledgment}
I would like to express my very great appreciation to professor Mehdi Golshani, for valuable discussions during my visit to the department of physics at IPM Tehran Iran 2010.

\end{document}